\documentclass{article}%
\usepackage{amsmath}
\usepackage{amsfonts}
\usepackage{amssymb}
\usepackage{graphicx}%
\newtheorem{theorem}{Theorem}

\newtheorem{corollary}[theorem]{Corollary}

\newtheorem{definition}[theorem]{Definition}
\newtheorem{example}[theorem]{Example}

\newtheorem{lemma}[theorem]{Lemma}

\newtheorem{proposition}[theorem]{Proposition}

\newtheorem{summary}[theorem]{Summary}
\newenvironment{proof}[1][Proof]{\noindent\textbf{#1.} }{\ \rule{0.5em}{0.5em}}
\begin{document}

\author{Diego E. Dominici\thanks{Department of Mathematics, Statistics and computer
Science, University of Illinois at Chicago (m/c 249), 851 South Morgan Street,
Chicago, IL 60607-7045, USA (ddomin1@uic.edu)}}
\title{The inverse of the cumulative standard normal probability function.}
\date{}
\maketitle

\begin{abstract}
Some properties of the inverse of the function $\ \ N(x)=\frac{1}{\sqrt{2\pi}%
}\int_{-\infty}^{x}e^{-\frac{t^{2}}{2}}dt$ \ are studied. Its derivatives,
integrals and asymptotic behavior are presented.

\end{abstract}

\section{Introduction}

It would be difficult to overestimate the importance of the standard normal
(or Gauss) distribution. It finds widespread application in almost every
scientific discipline, e.g., probability theory, the theory of errors, heat
conduction, biology, economics, physics, neural networks \cite{menon}, etc. It
plays a fundamental role in the financial mathematics, being part of the
Black-Scholes formula \cite{BS}, and its inverse is used in computing the
implied volatility of an option \cite{lee}. Yet, little is known about the
properties of the inverse function, e.g., series expansions, asymptotic
behavior, integral representations. The major work done has been in computing
fast and accurate algorithms for numerical calculations \cite{blair}.

Over the years a few articles have appeared with analytical studies of the
closely related error function%

\[
\operatorname{erf}(x)=\frac{2}{\sqrt{\pi}}\int\limits_{0}^{x}e^{-t^{2}}dt
\]
and its complement%

\[
\operatorname{erfc}(x)=\frac{2}{\sqrt{\pi}}\int\limits_{x}^{\infty}e^{-t^{2}%
}dt\text{ .}%
\]
Philip \cite{philip} introduced the notation \textquotedblleft inverfc($x$%
)\textquotedblright\ to denote the inverse of the complementary error
function. He gave the first terms in the power series for inverfc($x$),
asymptotic formulas for small $x$ in terms of continued logarithms, and some
expressions for the derivatives and integrals. Carlitz \cite{carlitz}, studied
the arithmetic properties of the coefficients in the power series of
inverfc($x$). Strecok \cite{strecok} computed the first 200 terms in the
series of inverfc($x$), and some expansions in series of Chebyshev
polynomials. Finally, Fettis \cite{fettis} studied inverfc($x$) for small $x$,
using an iterative sequence of logarithms.

The purpose of this paper is to present some new results on the derivatives,
integrals, and asymptotics of the inverse of the cumulative standard normal
probability function%

\[
N(x)=\frac{1}{\sqrt{2\pi}}\int_{-\infty}^{x}e^{-\frac{t^{2}}{2}}dt
\]
which we call $S(x)$. In section 2 we derive an ODE satisfied by $S(x)$, and
solve it using a power series. We introduce a family of polynomials $P_{n}$
related to the calculation of higher derivatives of $S(x)$. In section 3 we
study some properties of the $P_{n}$, such as relations between coefficients,
recurrences, and generating functions. We also derive a general formula for
$P_{n}$ using the idea of \textquotedblleft nested
derivatives\textquotedblright, and we compare the $P_{n}$ with the Hermite
polynomials $H_{n}$.

In section 4 we extend the definition of the $P_{n}$ to $n<0$ and use them to
calculate the integrals of $S(x)$. We also compute the integrals of powers of
$S(x)$ on the interval $[0,1]$. Section 5 is dedicated to asymptotics of
$S(x)$ for $x\rightarrow0,\ x\rightarrow1$ using the function Lambert W. With
the help of those formulas we derive an approximation to $S(x)$ valid in the
interval $[0,1]$ with error $\varepsilon$, $\left\vert \varepsilon\right\vert
<0.0023$. Finally, appendix A contains the first 20 non-zero coefficients in
the series of $S(x)$, and the first 10 polynomials $P_{n}$.

\section{Derivatives}

\begin{definition}
Let $S(x)$ denote the inverse of%
\[
N(x)=\frac{1}{\sqrt{2\pi}}\int_{-\infty}^{x}e^{-\frac{t^{2}}{2}}dt
\]

satisfying
\begin{equation}
S\circ N(x)=N\circ S(x)=x
\end{equation}

In terms of the error function $\operatorname{erf}(x)$,
\[
N(x)=\frac{1}{2}\left[  \operatorname{erf}\left(  \frac{x}{\sqrt{2}}\right)
+1\right]
\]

\end{definition}

\begin{proposition}
$S(x)$ satisfies the IVP%
\begin{gather}
S^{\prime\prime}=S(S^{\prime})^{2}\\
S\left(  \tfrac{1}{2}\right)  =0,\text{ }S^{\prime}\left(  \tfrac{1}%
{2}\right)  =\sqrt{2\pi}\nonumber
\end{gather}

\end{proposition}

\begin{proof}
Since $N(0)=1/2$, in follows that $S(1/2)=0$. From (2.1) we get%
\[
S^{\prime}[N(x)]=\frac{1}{N^{\prime}(x)}=\sqrt{2\pi}e^{\frac{x^{2}}{2}}%
=\sqrt{2\pi}e^{\frac{S^{2}[N(x)]}{2}}%
\]

Substituting $N(x)=y$ we have%
\[
S^{\prime}(y)=\sqrt{2\pi}e^{\frac{S^{2}(y)}{2}},\quad S^{\prime}(\frac{1}%
{2})=\sqrt{2\pi}%
\]

Differentiating $\ln[S^{\prime}(y)]$, we get (2.2).
\end{proof}

\begin{proposition}%
\[
S^{(n+2)}(x)=\sum_{i=0}^{n}\sum_{j=0}^{i}\binom{n}{i}\binom{i}{j}%
S^{(n-i)}(x)S^{(i-j+1)}(x)S^{(j+1)}(x),\quad n\geq0.
\]

\end{proposition}

\begin{proof}
Taking the n$^{\text{th}}$ derivative of (2.2), and using Leibnitz's Theorem,
we have%
\begin{align*}
S^{(n+2)}  &  =\sum_{i=0}^{n}\binom{n}{i}S^{(n-i)}(S^{\prime}S^{\prime}%
)^{(i)}\\
&  =\sum_{i=0}^{n}\binom{n}{i}S^{(n-i)}\sum_{j=0}^{i}\binom{i}{j}%
S^{(i-j+1)}S^{(j+1)}\\
&  =\sum_{i=0}^{n}\sum_{j=0}^{i}\binom{n}{i}\binom{i}{j}S^{(n-i)}%
S^{(i-j+1)}S^{(j+1)}.
\end{align*}

\end{proof}

\begin{corollary}
If $D_{n}=\frac{d^{n}S}{dx^{n}}(\frac{1}{2})$ , then%
\[
D_{2n}=0,\text{ \ }n\geq0.
\]

\end{corollary}

Putting $D_{n}=(2\pi)^{\frac{n}{2}}C_{n}$ , we can write%
\[
S(x)=\sum_{n\geq0}(2\pi)^{\frac{2n+1}{2}}\frac{C_{2n+1}}{(2n+1)!}(x-\frac
{1}{2})^{2n+1}%
\]

where%
\[
C_{1}=1,\text{ }C_{3}=1,\text{ }C_{5}=7,\text{ }C_{7}=127,\text{ }\ldots
\quad.
\]

\begin{proposition}%
\[
S^{(n)}=P_{n-1}(S)(S^{\prime})^{n}\quad n\geq1
\]

where $P_{n}(x)$ is a polynomial of degree $n$ satisfying the recurrence%
\begin{equation}
P_{0}(x)=1,\quad P_{n}(x)=P_{n-1}^{\prime}(x)+nxP_{n-1}(x),\quad n\geq1
\end{equation}

so that
\[
P_{1}(x)=x,\quad P_{2}(x)=1+2x^{2},\quad P_{3}(x)=7x+6x^{3},~\ldots~.
\]

\end{proposition}

\begin{proof}
We use induction on $n$. For $n$ = 2 the result follows from (2.2). If we
assume the result is true for $n$ then
\begin{align*}
S^{(n+1)}  &  =[P_{n-1}(S)(S^{\prime})^{n}]^{\prime}\\
&  =P_{n-1}^{\prime}(S)S^{\prime}(S^{\prime})^{n}+P_{n-1}(S)n(S^{\prime
})^{n-1}S^{\prime\prime}\\
&  =P_{n-1}^{\prime}(S)(S^{\prime})^{n+1}+P_{n-1}(S)n(S^{\prime}%
)^{n-1}S(S^{\prime})^{2}\\
&  =[P_{n-1}^{\prime}(S)+nSP_{n-1}(S)](S^{\prime})^{n+1}\\
&  =P_{n}(S)(S^{\prime})^{n+1}%
\end{align*}

Since \ $P_{n-1}(x)$ is a polynomial of degree $n-1$ by hypothesis, is clear
that%
\[
P_{n}(x)=P_{n-1}^{\prime}(x)+nxP_{n-1}(x)
\]

is a polynomial of degree $n$.
\end{proof}

\begin{corollary}%
\[
C_{n}=P_{n-1}(0)
\]

\end{corollary}

\section{The polynomials $P_{n}(x)$}

\begin{lemma}
If we write
\[
P_{n}(x)=\sum_{k=0}^{n}Q_{k}^{n}x^{k}%
\]

we have
\begin{align}
Q_{0}^{n}  &  =Q_{1}^{n-1}\nonumber\\
Q_{k}^{n}  &  =nQ_{k-1}^{n-1}+(k+1)Q_{k+1}^{n-1}\quad k=1,\ldots,n-2\\
Q_{k}^{n}  &  =nQ_{k}^{n-1}\quad k=n-1,n\quad\nonumber
\end{align}

\end{lemma}

\begin{proof}%
\begin{align*}
\sum_{k=0}^{n}Q_{k}^{n}x^{k}  &  =P_{n}\\
&  =\frac{d}{dx}P_{n-1}+nxP_{n-1}\\
&  =\sum_{k=0}^{n-1}Q_{k}^{n-1}kx^{k-1}+\sum_{k=0}^{n-1}nQ_{k}^{n-1}x^{k+1}\\
&  =\sum_{k=0}^{n-2}Q_{k+1}^{n-1}(k+1)x^{k}+\sum_{k=1}^{n}nQ_{k-1}^{n-1}x^{k}%
\end{align*}

\end{proof}

\begin{corollary}
In matrix form (3.1) reads $A^{(n)}Q^{n-1}=Q^{n}$ where $A^{(n)}$ $\in$
$\Re^{(n+1)\ \times\ n}$ is given by%
\[
A_{i,j}^{(n)}=\left\{
\begin{array}
[c]{c}%
i,\quad j=i+1,\quad i=1,\ldots,n-1\\
n,\quad j=i-1,\quad i=2,\ldots,n+1\\
0,\quad\quad ow\quad
\end{array}
\right.
\]

In other words, $A^{(n)}$ is a rectangular matrix with zeros in the diagonal,
the numbers $1,2,3,\ldots,n-1$ above it, the number $n$ below it and zeros
everywhere else%
\[
A^{(n)}=\left[
\begin{array}
[c]{cccccccc}%
0 & 1 & 0 & 0 & 0 & \cdots & 0 & 0\\
n & 0 & 2 & 0 & 0 & \cdots & 0 & 0\\
0 & n & 0 & 3 & 0 & \cdots & 0 & 0\\
\vdots & \vdots & \vdots & \vdots & \vdots & \ddots & \vdots & \vdots\\
0 & 0 & 0 & 0 & 0 & \cdots & 0 & n-1\\
0 & 0 & 0 & 0 & 0 & \cdots & n & 0\\
0 & 0 & 0 & 0 & 0 & \cdots & 0 & n
\end{array}
\right]
\]

and%
\[
Q^{n}=\left[
\begin{array}
[c]{c}%
Q_{0}^{n}\\
Q_{1}^{n}\\
Q_{2}^{n}\\
\vdots\\
Q_{n}^{n}%
\end{array}
\right]
\]

With the help of these matrices, we have an expression for the coefficients of
$P_{n}(x)$%
\[
Q^{n}=\prod\limits_{k=1}^{n}A^{(n-k+1)}=A^{(n)}A^{(n-1)}\cdots A^{(1)}%
\]

\end{corollary}

\begin{proposition}
The polynomials $P_{n}(x)$ satisfy the recurrence relation%
\[
P_{n+1}(x)=\sum\limits_{i=0}^{n}\sum\limits_{j=0}^{i}\binom{n}{i}\binom{i}%
{j}P_{n-i-1}(x)P_{i-j}(x)P_{j}(x)
\]

\end{proposition}

\begin{proof}%
\begin{align*}
&  P_{n+1}(S)(S^{\prime})^{n+2}\\
=  &  S^{(n+2)}=\sum_{i=0}^{n}\sum_{j=0}^{i}\binom{n}{i}\binom{i}{j}%
S^{(n-i)}S^{(i-j+1)}S^{(j+1)}\\
&  =\sum_{i=0}^{n}\sum_{j=0}^{i}\binom{n}{i}\binom{i}{j}P_{n-i-1}%
(S)(S^{\prime})^{n-i}P_{i-j}(S)(S^{\prime})^{i-j+1}P_{j}(S)(S^{\prime}%
)^{j+1}\\
&  =(S^{\prime})^{n+2}\sum_{i=0}^{n}\sum_{j=0}^{i}\binom{n}{i}\binom{i}%
{j}P_{n-i-1}(S)P_{i-j}(S)P_{j}(S).
\end{align*}

\end{proof}

\begin{proposition}
The exponential generating function of the polynomials $P_{n}(x)$ is%
\[
e^{\frac{1}{2}S^{2}[N(x)+tN^{\prime}(x)]-\frac{x^{2}}{2}}=F(x,t)=\sum_{k\geq
0}P_{k}(x)\frac{t^{k}}{k!}%
\]

\end{proposition}

\begin{proof}
Since
\begin{align*}
F(x,t)  &  =\sum_{k\geq0}P_{k}(x)\frac{t^{k}}{k!}\\
&  =1+\sum_{k\geq1}\frac{d}{dx}P_{k-1}(x)\frac{t^{k}}{k!}+\sum_{k\geq
1}kxP_{k-1}(x)\frac{t^{k}}{k!}\\
&  =1+\sum_{k\geq0}\frac{d}{dx}P_{k}(x)\frac{t^{k+1}}{(k+1)k!}+\sum_{k\geq
0}xP_{k}(x)\frac{t^{k+1}}{k!}\\
&  =1+\int_{0}^{t}\frac{\partial}{\partial x}F(x,s)ds+xtF(x,t)
\end{align*}

it follows that $F(x,t)$ satisfies the differential-integral equation%
\begin{equation}
1+(xt-1)F(x,t)+\frac{\partial}{\partial x}\int_{0}^{t}F(x,s)ds=0
\end{equation}

Differentiating (3.2) with respect to $t$ we get%
\[
xF(x,t)+(xt-1)\frac{\partial}{\partial t}F(x,t)+\frac{\partial}{\partial
x}F(x,t)=0
\]

whose general solution is of the form%
\[
F(x,t)=e^{-\frac{x^{2}}{2}}G\left(  te^{-\frac{x^{2}}{2}}+\sqrt{2\pi}\left[
N(x)-\frac{1}{2}\right]  \right)
\]

for some function $G(z)$.

From (3.2) we know that $F(x,0)=1$, and hence%
\[
G\left(  \sqrt{2\pi}\left[  N(x)-\frac{1}{2}\right]  \right)  =e^{\frac{x^{2}%
}{2}},
\]

which implies that%
\[
G(z)=e^{\frac{1}{2}S^{2}\left(  \frac{z}{\sqrt{2\pi}}+\frac{1}{2}\right)  }.
\]

Therefore,%
\[
F(x,t)=e^{\frac{1}{2}S^{2}[N(x)+tN^{\prime}(x)]-\frac{x^{2}}{2}}%
\]

\end{proof}

\begin{definition}
We define the \textquotedblleft nested derivative\textquotedblright%
\ $\mathfrak{D}^{(n)}$\ by%
\begin{align*}
\mathfrak{D}^{(0)}[f](x)  &  \equiv1\\
\mathfrak{D}^{(n)}[f]\,(x)  &  =\frac{d}{dx}\left\{  f(x)\times\mathfrak{D}%
^{(n-1)}[f]\,(x)\right\}  ,\quad n\geq1
\end{align*}

\end{definition}

\begin{example}
\begin{enumerate}
\item
\[
\mathfrak{D}^{(n)}\left[  e^{ax}\right]  =n!a^{n}e^{nax}%
\]

\item
\[
\mathfrak{D}^{(n)}\left[  x\right]  =1
\]

\item
\[
\mathfrak{D}^{(n)}\left[  x^{2}\right]  =(n+1)!x^{n}%
\]

\end{enumerate}
\end{example}

\begin{proposition}%
\[
P_{n}(x)=e^{-\frac{n}{2}x^{2}}\mathfrak{D}^{(n)}\left[  e^{\frac{x^{2}}{2}%
}\right]
\]

\end{proposition}

\begin{proof}
We use induction on $n$. For $n=0$ the result follows from the definition of
$\mathfrak{D}^{(n)}$. Assuming the result is true for $n-1$
\begin{align*}
P_{n}(x)  &  =P_{n-1}^{\prime}(x)+nxP_{n-1}(x)\\
&  =\frac{d}{dx}[e^{-\frac{(n-1)}{2}x^{2}}\mathfrak{D}^{(n-1)}(e^{\frac{x^{2}%
}{2}})]+nxe^{-\frac{(n-1)}{2}x^{2}}\mathfrak{D}^{(n-1)}(e^{\frac{x^{2}}{2}})\\
&  =-(n-1)xe^{-\frac{(n-1)}{2}x^{2}}\mathfrak{D}^{(n-1)}(e^{\frac{x^{2}}{2}%
})+e^{-\frac{(n-1)}{2}x^{2}}\frac{d}{dx}[\mathfrak{D}^{(n)}(e^{\frac{x^{2}}%
{2}})]+\\
&  +nxe^{-\frac{(n-1)}{2}x^{2}}\mathfrak{D}^{(n-1)}(e^{\frac{x^{2}}{2}})\\
&  =e^{-\frac{(n-1)}{2}x^{2}}\left[  x\mathfrak{D}^{(n-1)}(e^{\frac{x^{2}}{2}%
})+\frac{d}{dx}\mathfrak{D}^{(n-1)}(e^{\frac{x^{2}}{2}})\right] \\
&  =e^{-\frac{(n-1)}{2}x^{2}}e^{-\frac{1}{2}x^{2}}\frac{d}{dx}\left[
e^{\frac{1}{2}x^{2}}\mathfrak{D}^{(n-1)}(e^{\frac{x^{2}}{2}})\right] \\
&  =e^{-\frac{n}{2}x^{2}}\mathfrak{D}^{(n)}\left[  e^{\frac{x^{2}}{2}}\right]
\end{align*}

\end{proof}

\begin{summary}
We conclude this section by comparing the properties of $P_{n}(x)$ with the
well known formulas for the Hermite polynomials $H_{n}(x)$ \cite{lebedev}.
Since the $H_{n}$ are deeply related with the function $N(x)$, we would expect
to see some similarities between the $H_{n}$ and the $P_{n}$.%

\[%
\begin{tabular}
[c]{|c|c|}\hline
$P_{n}(x)$ & $H_{n}(x)$\\
& \\\hline
\multicolumn{1}{|l|}{$P_{n}(x)=e^{-\frac{n}{2}x^{2}}\mathfrak{D}%
^{(n)}(e^{\frac{x^{2}}{2}})$} & \multicolumn{1}{|l|}{$H_{n}(x)=(-1)^{n}%
e^{x^{2}}\frac{d^{n}}{dx^{n}}(e^{-x^{2}})$}\\
\multicolumn{1}{|l|}{} & \multicolumn{1}{|l|}{}\\\hline
\multicolumn{1}{|l|}{$\sum\limits_{k\geq0}P_{k}(x)\frac{t^{k}}{k!}=e^{\frac
{1}{2}S^{2}[N(x)+tN^{\prime}(x)]-\frac{x^{2}}{2}}$} &
\multicolumn{1}{|l|}{$\sum\limits_{k\geq0}H_{k}(x)\frac{t^{k}}{k!}%
=e^{2xt-t^{2}}$}\\
\multicolumn{1}{|l|}{} & \multicolumn{1}{|l|}{}\\\hline
\multicolumn{1}{|l|}{$P_{n}(x)=P_{n-1}^{\prime}(x)+nxP_{n-1}(x)$} &
\multicolumn{1}{|l|}{$H_{n}(x)=-H_{n-1}^{\prime}(x)+2xH_{n-1}(x)$}\\
\multicolumn{1}{|l|}{} & \multicolumn{1}{|l|}{}\\\hline
\multicolumn{1}{|l|}{$S^{(n)}=P_{n-1}(S)(S^{\prime})^{n}$} &
\multicolumn{1}{|l|}{$N^{(n)}=\left(  \frac{-1}{\sqrt{2}}\right)
^{n-1}H_{n-1}\left(  \frac{x}{\sqrt{2}}\right)  N^{\prime}$}\\
\multicolumn{1}{|l|}{} & \multicolumn{1}{|l|}{}\\\hline
\end{tabular}
\ \ \
\]

\end{summary}

\section{Integrals of S(x)}

\begin{definition}%
\[
S^{(-n)}(x)=\int\limits_{0}^{x}\int\limits_{0}^{x_{1}}\cdots\int
\limits_{0}^{x_{n-1}}S(x_{n})\ dx_{n}dx_{n-1}\ldots dx_{1},\qquad n\geq1.
\]

\end{definition}

\begin{lemma}%
\begin{equation}
P_{n-1}(x)=e^{-\frac{n}{2}x^{2}}\left[  P_{n-1}(0)+\int\limits_{0}^{x}%
e^{\frac{n}{2}t^{2}}P_{n}(t)dt\right]
\end{equation}

\end{lemma}

\begin{proof}
It follows immediately from solving the ODE for $P_{n-1}$ in terms of $P_{n}$.
\end{proof}

\begin{proposition}
Using (4.1 ) to define $P_{n}(x)$ for $n<0$ yields%
\begin{align}
P_{-1}(x)  &  =x\nonumber\\
P_{-2}(x)  &  =-1\\
P_{-3}(x)  &  =-\sqrt{\pi}e^{x^{2}}N\left(  \sqrt{2}x\right) \nonumber
\end{align}

and the relation%
\[
S^{(n)}=P_{n-1}(S)(S^{\prime})^{n}%
\]

still holds.
\end{proposition}

\begin{proof}
For $n=0$ we have
\begin{align*}
P_{-1}(x)  &  =P_{-1}(0)+x\\
S  &  =S^{(0)}=P_{-1}(S)
\end{align*}

so $P_{-1}(0)=0$. For $n=-1$%
\[
P_{-2}(x)=e^{\frac{x^{2}}{2}}\left[  P_{-2}(0)+1\right]  -1
\]

We can calculate $S^{(-1)}$ explicitly by%
\begin{align*}
S^{(-1)}(x)  &  =\int\limits_{0}^{x}S(t)dt=\frac{1}{\sqrt{2\pi}}%
\int\limits_{-\infty}^{S(x)}S[N(z)]e^{-\frac{z^{2}}{2}}dz\\
&  =\frac{1}{\sqrt{2\pi}}\int\limits_{-\infty}^{S(x)}ze^{-\frac{z^{2}}{2}}dz\\
&  =-\frac{1}{\sqrt{2\pi}}e^{-\frac{S(x)^{2}}{2}}=-(S^{\prime})^{-1}%
\end{align*}

Hence, $P_{-2}(0)=-1$.

Finally, for $n=-2$%
\[
P_{-3}(x)=e^{x^{2}}\left[  P_{-3}(0)-\sqrt{\pi}N(\sqrt{2}x)+\frac{\sqrt{\pi}%
}{2}\right]
\]

A similar calculation as the one above, making a change of variables $t=N(z)$
in the integral of $S^{(-1)}(x)$ yields%
\[
S^{(-2)}(x)=-\frac{1}{2\sqrt{\pi}}N[\sqrt{2}S(x)]
\]

and we conclude that $P_{-3}(0)=-\frac{\sqrt{\pi}}{2}$.
\end{proof}

\begin{corollary}%
\[
S\left[  -2\sqrt{\pi}S^{(-2)}\right]  =\sqrt{2}S
\]

\end{corollary}

\begin{proposition}%
\[
\int\limits_{0}^{1}S^{n}(x)dx=\left\{
\begin{array}
[c]{c}%
\prod\limits_{i=1}^{k}(2i+1),\quad n=2k,\quad k\geq1\\
0,\quad n=2k+1,\quad k\geq0
\end{array}
\right.
\]

\end{proposition}

\begin{proof}%
\begin{align*}
\int\limits_{0}^{1}S^{n}(x)dx  &  =\int\limits_{-\infty}^{\infty}%
z^{n}N^{\prime}(z)dz\\
&  =\frac{1}{\sqrt{2\pi}}\int\limits_{-\infty}^{\infty}z^{n}e^{-\frac{1}%
{2}z^{2}}dz\\
&  =\prod\limits_{i=1}^{k}(2i+1),\quad n=2k,\quad k\geq1
\end{align*}

\end{proof}

\section{Asymptotics}

\begin{definition}
We'll denote by $LW(x)$ the function Lambert W \cite{corless},%
\begin{equation}
LW(x)e^{LW(x)}=x
\end{equation}

This function has the series representation \cite{corless1}%
\[
LW(x)=\sum_{n\geq1}\frac{(-n)^{n-1}}{n!}x^{n},
\]

the derivative%
\[
\frac{d}{dx}LW=\frac{LW(x)}{x[1+LW(x)]}\quad\text{if }x\neq0,
\]

and it has the asymptotic behavior%
\[
LW(x)\sim\ln(x)-\ln[\ln(x)]\quad x\rightarrow\infty.
\]

\end{definition}

\begin{proposition}%
\begin{align*}
S(x)  &  \sim g_{0}(x)=-\sqrt{LW\left(  \frac{1}{2\pi x^{2}}\right)  },\quad
x\rightarrow0\\
S(x)  &  \sim g_{1}(x)=\sqrt{LW\left(  \frac{1}{2\pi(x-1)^{2}}\right)  },\quad
x\rightarrow1
\end{align*}

Both functions $g_{0}(x)$ and $g_{1}(x)$ satisfy the ODE%
\[
g^{\prime\prime}=g(g^{\prime})^{2}\left[  1+\frac{2}{g^{2}(1+g^{2})}\right]
\sim g(g^{\prime})^{2},\quad\text{for \ }\left\vert g\right\vert
\rightarrow\infty
\]

\end{proposition}

\begin{proof}%
\begin{align*}
N(x)  &  \sim\frac{1}{\sqrt{2\pi}}e^{-\frac{x^{2}}{2}}\frac{1}{x},\quad
x\rightarrow-\infty\\
t  &  \sim\frac{1}{\sqrt{2\pi}}e^{-\frac{S(t)^{2}}{2}}\frac{1}{S(t)},\quad
t\rightarrow0\\
S(t)e^{\frac{S(t)^{2}}{2}}  &  \sim\frac{1}{\sqrt{2\pi}t},\quad t\rightarrow
0\\
S^{2}(t)e^{S^{2}(t)}  &  \sim\frac{1}{2\pi t^{2}},\quad t\rightarrow0
\end{align*}

Using the definition of $LW(x)$ we have
\[
S^{2}(t)\sim LW\left(  \frac{1}{2\pi t^{2}}\right)  ,\quad t\rightarrow0
\]

or
\[
S(t)\sim-\sqrt{LW\left(  \frac{1}{2\pi t^{2}}\right)  },\quad t\rightarrow0
\]

The case $x\rightarrow1$ is completely analogous.
\end{proof}

\begin{corollary}
Combining the above expressions, we can get the approximation%
\begin{equation}
S(x)\simeq g_{2}(x)=(2x-1)\sqrt{LW\left(  \frac{1}{2\pi x^{2}(x-1)^{2}%
}\right)  }%
\end{equation}

good through the interval $(0,1)$.

We can refine it even more by putting%
\begin{align}
S(x)  &  \simeq g_{3}(x)=Q(x)\sqrt{LW\left(  \frac{1}{2\pi x^{2}(x-1)^{2}%
}\right)  }\\
Q(x)  &  =-1+(6-2\sqrt{2\pi})x+(-12+6\sqrt{2\pi})x^{2}+(8-4\sqrt{2\pi}%
)x^{3}\nonumber
\end{align}

where $Q(x)$ has been chosen such that%
\[
Q(0)=-1,~Q(1)=0,~Q(1/2)=0,~Q^{\prime}(1/2)=\sqrt{2\pi},~Q^{\prime\prime
}(1/2)=0
\]

\end{corollary}

\section{Appendix}

The first 10 $P_{n}(x)$ are%
\begin{align*}
P_{0}(x)  &  =1\\
P_{1}(x)  &  =x\\
P_{2}(x)  &  =1+2x^{2}\\
P_{3}(x)  &  =7x+6x^{3}\\
P_{4}(x)  &  =7+46x^{2}+24x^{4}\\
P_{5}(x)  &  =127x+326x^{3}+120x^{5}\\
P_{6}(x)  &  =127+1740x^{2}+2556x^{4}+720x^{6}\\
P_{7}(x)  &  =4369x+22404x^{3}+22212x^{5}+5040x^{7}\\
P_{8}(x)  &  =4369+102164x^{2}+290292x^{4}+212976x^{6}+40320x^{8}\\
P_{9}(x)  &  =243649x+2080644x^{3}+3890484x^{5}+2239344x^{7}+362880x^{9}\\
P_{10}(x)  &  =243649+8678422x^{2}+40258860x^{4}+54580248x^{6}+\\
&  +25659360x^{8}+3628800x^{10}%
\end{align*}

The first few odd $C_{n}$ are%

\begin{tabular}
[c]{cc}%
$n$ & $C_{n}$\\
1 & 1\\
3 & 1\\
5 & 7\\
7 & 127\\
9 & 4369\\
11 & 243649\\
13 & 20036983\\
15 & 2280356863\\
17 & 343141433761\\
19 & 65967241200001\\
21 & 15773461423793767\\
23 & 4591227123230945407\\
25 & 1598351733247609852849\\
27 & 655782249799531714375489\\
29 & 313160404864973852338669783\\
31 & 172201668512657346455126457343\\
33 & 108026349476762041127839800617281\\
35 & 76683701969726780307420968904733441\\
37 & 61154674195324330125295778531172438727\\
39 & 54441029530574028687402753586278549396607\\
41 & 53789884101606550209324949796685518122943569
\end{tabular}


\begin{thebibliography}{99}                                                                                               %


\bibitem {blair}J.~M. Blair, C.~A. Edwards, and J.~H. Johnson. \newblock
Rational {C}hebyshev approximations for the inverse of the error function.
\newblock {\em Math. Comp.}, 30(136):827--830, 1976.

\bibitem {BS}Fisher Black and Myron~S Scholes. \newblock The pricing of
options and corporate liabilities.
\newblock {\em Journal of Political Economy}, 81(3):637--654, 1973.

\bibitem {carlitz}L.~Carlitz. \newblock The inverse of the error function.
\newblock {\em Pacific J. Math.}, 13:459--470, 1963.

\bibitem {corless}R.~M. Corless, G.~H. Gonnet, D.~E.~G. Hare, D.~J. Jeffrey,
and D.~E. Knuth. \newblock On the {L}ambert ${W}$ function.
\newblock {\em Adv. Comput. Math.}, 5(4):329--359, 1996.

\bibitem {corless1}Robert~M. Corless, David~J. Jeffrey, and Donald~E. Knuth.
\newblock A sequence of series for the {L}ambert ${W}$ function. \newblock In
\emph{Proceedings of the 1997 International Symposium on Symbolic and
Algebraic Computation (Kihei, HI)}, pages 197--204 (electronic), New York,
1997. ACM.

\bibitem {fettis}Henry~E. Fettis. \newblock A stable algorithm for computing
the inverse error function in the \textquotedblleft tail-end\textquotedblright%
\ region. \newblock {\em Math. Comp.}, 28:585--587, 1974.

\bibitem {lebedev}N.~N. Lebedev.
\newblock {\em Special functions and their applications}. \newblock
Prentice-Hall Inc., Englewood Cliffs, N.J., 1965.

\bibitem {lee}C.F. Lee and A.~Tucker. \newblock An alternative method for
obtaining the implied standard deviation.
\newblock {\em The Journal of Financial Engineering}, 1:369--375, 1992.

\bibitem {menon}A.~Menon, K.~Mehrotra, C.~K. Mohan, and S.~Ranka. \newblock
Characterization of a class of sigmoid functions with applications to neural
networks. \newblock {\em Neural Networks}, 9(5):819--835, 1996.

\bibitem {philip}J.~R. Philip. \newblock The function inverfc $\theta$.
\newblock {\em Austral. J. Phys.}, 13:13--20, 1960.

\bibitem {strecok}Anthony Strecok. \newblock On the calculation of the inverse
of the error function. \newblock {\em Math. Comp.}, 22:144--158, 1968.
\end{thebibliography}
\end{document}